\documentclass[12pt]{article}
\usepackage{natbib}
\usepackage{graphicx,
setspace,amsmath,amssymb, mathrsfs,subfigure,url,color}
\newenvironment{frcseries}{\fontfamily{frc}\selectfont}{}

\usepackage{epstopdf}
\newcommand{\ep}{\epsilon}

\newcommand{\bay}{\begin{array}}
\newcommand{\eay}{\end{array}}

\newcommand{\bqa}{\begin{eqnarray*}}
\newcommand{\eqa}{\end{eqnarray*}}
\newcommand{\bqan}{\begin{eqnarray}}
\newcommand{\eqan}{\end{eqnarray}}
\newcommand{\bqt}{\begin{quote}}
\newcommand{\eqt}{\end{quote}}
\newcommand{\bt}{\begin{tabbing}}
\newcommand{\et}{\end{tabbing}}
\newcommand{\bit}{\begin{itemize}}
\newcommand{\eit}{\end{itemize}}
\newcommand{\ben}{\begin{enumerate}}
\newcommand{\een}{\end{enumerate}}
\newcommand{\beq}{\begin{equation}}
\newcommand{\eeq}{\end{equation}}
\newcommand{\bdefi}{\begin{definition}}
\newcommand{\edefi}{\end{definition}}
\newcommand{\bpro}{\begin{proposition}}
\newcommand{\epro}{\end{proposition}}
\newcommand{\bco}{\begin{corollary}}
\newcommand{\eco}{\end{corollary}}
\newcommand{\bdes}{\begin{description}}
\newcommand{\edes}{\end{description}}

\def\wt{\widetilde}

\def\wh{\widehat}

\def\log{\hbox{log}}

\def\boxit#1{\vbox{\hrule\hbox{\vrule\kern6pt
          \vbox{\kern6pt#1\kern6pt}\kern6pt\vrule}\hrule}}

\def\bse{\begin{eqnarray*}}
\def\ese{\end{eqnarray*}}
\def\be{\begin{eqnarray}}
\def\ee{\end{eqnarray}}
\def\bq{\begin{equation}}
\def\eq{\end{equation}}

\def\wh{\widehat}

\def\tr{\textcolor{black}}

\newtheorem{proposition}{Proposition}

\newcommand{\blem}{\begin{lemma}}
\newcommand{\elem}{\end{lemma}}
\newcommand{\bthe}{\begin{theorem}}
\newcommand{\ethe}{\end{theorem}}

\newtheorem{definition}{Definition}[section]
\newtheorem{lemma}[definition]{Lemma}
\newtheorem{theorem}[definition]{Theorem}

\def\delete#1{\iffalse #1 \fi}

\def\bse{\begin{eqnarray*}}
\def\ese{\end{eqnarray*}}
\def\bee{\begin{enumerate}}
\def\eee{\end{enumerate}}
\def\bqe{\begin{eqnarray}}
\def\eqe{\end{eqnarray}}
\def\bed{\begin{description}}
\def\eed{\end{description}}
\def\bei{\begin{itemize}}
\def\eei{\end{itemize}}

\def\pmb#1{\setbox0=\hbox{#1}%
    \kern-.025em\copy0\kern-\wd0
    \kern.05em\copy0\kern-\wd0
    \kern-.025em\raise.0433em\box0 }
\def\pmbh#1#2{\setbox0=\hbox{#1}%
    \setbox1=\hbox{#2}%
    \kern-.025em\copy0\kern-\wd0
    \kern.05em\copy1\kern-\wd0
    \kern-.025em\raise.0433em\box0 }

\def\frac#1#2{{#1\over#2}}

\def\boxit#1{\vbox{\hrule\hbox{\vrule\kern6pt
   \vbox{\kern6pt#1\kern6pt}\kern6pt\vrule}\hrule}}

\def\listing#1{\vskip 4mm\begin{verbatim}\input#1 \vskip 4mm}
\def\thick#1{\hbox{\rlap{$#1$}\kern0.25pt\rlap{$#1$}\kern0.25pt$#1$}}

\def\wt{\widetilde}
\def\wh{\widehat}

\def\tr{{\mbox{\rm tr}}}

\def\pmbh{{\pmb h}}

\def\calE{{\cal E}}

\def\calH{{\cal H}}

\def\calL{{\cal L}}

\def\calN{{\cal N}}

\def\calT{{\cal T}}

\renewcommand\today{\ifcase\month\or
   Jan\or Feb\or Mar\or Apr\or May\or
   Jun\or Jul\or Aug\or Sep\or Oct\or Nov\or
   Dec\fi
   \space\number\day, \number\year}

\renewcommand{\log}{{\rm log}}

\newtheorem{thm}{Theorem}
\newtheorem{lem}{Lemma}

\newtheorem{rmk}{Remark}

\DeclareMathOperator*{\argmin}{arg\,min}
\begin{document}

\title{Functional Adaptive Huber Linear Regression}        

\author{
Ling Peng$^{a,b,c}$, Xiaohui Liu$^{b,a}$ Heng Lian$^{d,e}$\footnote{henglian@cityu.edu.hk}
\\
\begin{tabular}{l} 
{\small\it$^a$ School of Statistics and Data Science, Jiangxi University of Finance and Economics}\\
{\small\it $\;\; $ Nanchang, Jiangxi, China}\\
{\small\it$^b$ Key Laboratory of Data Science in Finance and Economics}\\
{\small\it $\;\; $  Jiangxi University of Finance and Economics, Nanchang, Jiangxi, China}\\
{\small\it$^c$ School of Mathematics and Statistics, Victoria University of Wellington}\\
{\small\it $\;\; $ Wellington 6140, New Zealand}\\
{\small\it$^d$ City University of Hong Kong Shenzhen Research Institute, Shenzhen, China}\\
{\small\it$^e$ Department of Mathematics, City University of Hong Kong, Hong Kong, China}
\end{tabular}
} 
\date{}          
\maketitle
\begin{abstract} Robust estimation has played an important role in statistical and machine learning. However, its applications to functional  linear regression are still under-developed. In this paper, we focus on Huber's loss with a diverging robustness parameter which was previously used in parametric models. Compared to other robust methods such as median regression, the distinction is that the proposed method aims to estimate the conditional mean robustly, instead of estimating the conditional median. We only require $(1+\kappa)$-th moment assumption ($\kappa>0$) on the noise distribution, and the established error bounds match the optimal rate in the least-squares case as soon as $\kappa\ge 1$. We establish convergence rate in probability when the functional predictor has a finite 4-th moment, and finite-sample bound with exponential tail  when the functional predictor is Gaussian, in terms of both prediction error and $L^2$ error. The results also extend to the case of functional estimation in a reproducing kernel Hilbert space (RKHS).

\end{abstract}

\section{Introduction}
Functional linear regression (FLR) is a powerful statistical framework that extends traditional linear regression models to accommodate functional predictors, allowing for a more flexible and nuanced analysis of complex data structures. The integration of functional data, often represented as curves, offers a comprehensive perspective on the underlying relationships between variables, making FLR a valuable tool in various scientific disciplines, including economics, biology, and signal processing.

Technically, there are many existing estimation approaches for functional linear regression \citep{yao05,hall07,crambes09}. It can be based on expanding the functional object in either a fixed basis or an estimated basis (such as using functional PCA). Alternatively, it can be fitted by using Tikhonov regularization. Although all these procedures could be robustified, we focus on the Tikhonov regularization approach, which works for function coefficient either in an $L^2$ space or in an RKHS space, with a unified proof.

Despite its merits, FLR faces challenges when confronted with outliers or heavy-tailed errors, which can adversely impact the estimation of model parameters and compromise the reliability of predictions. Traditional linear regression methods, based on the assumption of normally distributed errors, are particularly sensitive to deviations from this assumption. Consequently, the robustness of FLR becomes a critical consideration, especially when dealing with  data that may exhibit non-Gaussian and heavy-tailed characteristics. 

This paper proposes a new approach in the FLR setting to enhance the robustness through the incorporation of an adaptive Huber loss function. The Huber loss, renowned for its ability to balance the advantages of mean squared error and mean absolute error, is well-suited for situations where the data distribution may deviate from normality. By adapting the Huber loss function within the FLR framework, we aim to develop a method that is more resilient to the impact of outliers and capable of providing robust estimates of functional regression parameters. Traditionally, the robustness parameter in the Huber's loss (the point where the quadratic loss transitions to the absolute deviation loss) is usually treated as a fixed constant in the theoretical analysis, causing a nonignorable bias in conditional mean estimation.   Adaptive Huber's loss has been investigated for parametric mean regression models in \cite{adahuber2020}, in which the authors showed that by using a diverging robustness parameter, only $(1+\kappa)$-th moment with $\kappa>0$ is required for consistent estimation, while if $\kappa\ge 1$ the estimator is as efficient as the least squares estimator. 
Our aim is to investigate theoretically the performances of adaptive Huber regression for functional linear models. The proof in the functional setting is drastically different from that in the parametric setting, and it is thus very challenging to establish the optimal rate. In particular, the bias-variance tradeoff in the infinite-dimensional functional setting behaves very differently from that in the parametric model, and it is thus unclear whether $\kappa=1$ is still sufficient to achieve the same rate as for least-squares-based FLR.


We now discuss some related theoretical works on functional linear regression, also emphasizing our contribution.
\begin{itemize}
\item (Works on FLR) Theoretical properties of FLR have been investigated in numerous works, including but not limited to \cite{hall07,yuancai10,caiyuan12,mas13}. In these works, assumptions on finite second moment of the noise term are used to show (optimal) rates in expectation or in probability. On the other hand, although we did not directly find any such results in the literature of FLR, rates for the estimator with an exponential tail bound most likely require an exponential tail bound for the error (for example, when the error is sub-Gaussian). For the closely related problem of kernel ridge regression, \cite{caponnetto07} established exponential tail bound under such an assumption on the error distribution, and \cite{Steinwart2009} used an even stronger assumption that the error is bounded. In contrast, we can establish the optimal statistical rate in probability when $\|X\|$ is assumed to have 4-th bounded moments, while when $X$ is Gaussian we can establish bounds with exponential tail probability (without exponential tail assumptions on \emph{errors}). 

\item (Works on quantile FLR) Although quantile estimation is robust against heavy-tailed errors, and functional quantile regression has been studied in several works \citep{Kato2012,Chen2012,Li2021}, the main conceptual difference is that it focuses on the conditional median rather than the conditional mean. When the error distribution is asymmetric and/or heterogeneous, the two quantities can be quite different and thus quantile regression is not consistent for mean estimation.

\item (Work on Huber FLR) In the paper \cite{tong2023}, the author studied exactly the same problem in the RKHS setting. However, the rate obtained is not optimal in the sense that for any finite $\kappa$ (in the moment assumption of the error), the rate is slower than that of least square regression, and the two become the same only when $\kappa\rightarrow\infty$. We believe this is due to that the proof of \cite{tong2023} relies on the fact that Huber's loss $H_\gamma$ (as defined in the next section) is Lipschitz with Lipschitz constant $\gamma$, when applying the contraction inequality for Rademacher processes. This is a standard proof strategy for dealing with Lipschitz losses (for example used also in \cite{liannips2022} for the quantile loss). However, as $\gamma$ is diverging, this would introduce an extra $\gamma$ factor in their bound, leading to a rate that is sub-optimal. Furthermore, \cite{tong2023} only considers the case $r=0$ ($r$ as defined in our Assumption (B) is a smoothness parameter for the true functional parameter). Our proof does not use the Lipschitz property of the loss in this way, and is thus quite different from the proof strategy used in previous works on Lipschitz losses \citep{geoffrey2020,tong2023}.
\end{itemize}

The rest of the article is organized as follows.
In Section \ref{sec:huber}, we present the model based on adaptive Huber's loss and the assumptions and main theoretical results. Section \ref{sec:rkhs} briefly considers the RKHS setting for functional regression,  and we conclude in Section \ref{sec:con}.  

\textit{Notations}. We use $P_n$ to denote the empirical measure, and $P_n-P$ the centered empirical measure. We will consider functions in the Hilbert space $L^2(\calT)$ (we focus on $\calT=[0,1]$ in our paper, but it could be a more general compact metric space). The standard inner product and norm on $L^2(\calT)$ are denoted by $\langle.,.\rangle$ and $\|.\|$, respectively. For an operator on $L^2(\calT)$, $\|.\|_{HS}$, $\|.\|_{op}$ and ${\rm tr}(.)$ denote the Hilbert-Schmidt norm, the spectral norm and the trace norm, respectively. Recall that the Hilbert-Schmidt norm of an operator $A: L^2(\calT)\rightarrow L^2(\calT)$ is, given any orthonormal basis for $L^2(\calT)$  denoted by $\{e_1,e_2,\ldots\}$, 
$
\|A\|_{HS}=\sqrt{\sum_{j,k=1}^\infty \langle Ae_j,e_k\rangle^2},
$
and the trace norm is 
$
{\rm tr}(A)=\sum_{j=1}^\infty \langle Ae_j,e_j\rangle.
$
$A$ is a Hilbert-Schmidt operator (trace-class operator) if its Hilbert-Schmidt norm (trace norm) is finite. For a sequence of random variables $X_n$, we say $X_n$ converges (to zero) in probability with rate $r_n$, if for any $\delta\in(0,1)$, there exists some $C>0$ such that $P(|X_n|>Cr_n)<\delta$, $\forall n$. We say $X_n$ converges (to zero) in probability with exponential tail and rate $r_n$, if there exists some $C>0, p>0$ such that for any $\delta\in(0,1)$, we have
$P(|X_n|>Cr_n\log^p(\frac{1}{\delta}))<\delta$. Note that last equation can also be written as $P(|X_n|>Cur_n)<e^{-u^{1/p}}$, $\forall u>0$.

\section{Functional linear regression based on adaptive Huber's loss}\label{sec:huber}
The classical functional linear regression model imposes that, for an independent and identically distributed (i.i.d.) sample $(X_i,y_i)$, $i=1,\ldots,n$, 
\be\label{eqn:flr}
y_i=a+\langle f,X_i\rangle+\epsilon_i,
\ee
where $X_i$ is a random element in some Hilbert space containing functions defined on a compact metric space $\calT$, $f$ is the functional parameter to be estimated, $a\in\mathbb{R}$ is the intercept parameter, and $\epsilon_i$ is the mean-zero noise. Without loss of much generality, following the literature of FLR, we always assume $\calT=[0,1]$ equipped with the Lebesgue measure, and the Hilbert space is $L^2(\calT)$ and thus $\langle f,X_i\rangle=\int_\calT X_i(t)f(t)dt$.

With a regularization parameter $\lambda>0$, the standard estimator based on Tikhonov regularization is the minimizer of 
\bse
\min_{f\in L^2(\calT)}\frac{1}{n}\sum_{i=1}^n  (y_i-\int_\calT X_i(t)f(t)dt)^2+\lambda\|f\|^2.
\ese
As a robust version, we consider the estimator
\bse
\wh f=\argmin_{f\in L^2(\calT)}\frac{1}{n}\sum_{i=1}^n H_\gamma(y_i-\int_\calT X_i(t)f(t)dt)+\lambda\|f\|^2,
\ese
where 
\bse
H_\gamma(u)=\left\{\begin{array}{cc}
 \frac{u^2}{2} & |u|\le \gamma\\
 \gamma|u|-\frac{\gamma^2}{2} & |u|>\gamma.
 \end{array}\right.
\ese
We use $H_\gamma'$ to denote the first derivative of $H_\gamma$, and since we also want to regard the loss as a function of $f$, we denote $\nabla H_\gamma(f)=-H_\gamma'(y-\langle X,f\rangle)X$, which is actually the Fr\'echet derivative of $H_\gamma(y-\langle X,f\rangle)$ with respect to $f$.

In the definition of Huber's loss, $\gamma$ is called the robustness parameter. Note that $\gamma$ is the transition point between a quadratic loss function and an absolute deviation loss function. When $\gamma\rightarrow \infty$, it approaches the quadratic loss and when $\gamma\rightarrow 0$, it approaches the absolute deviation loss. A larger value of $\gamma$ makes it less biased for estimating the conditional mean, while being less robust. We will see that an appropriately diverging $\gamma$ as $n\rightarrow\infty$ will provide the optimal trade-off between the two. 

We will assume $E[X]=0$ and $E[y]=a=0$ for simplicity of notation. Our results can be extended to the general case since the intercept and $E[X]$ can be estimated easily with faster rates. Let $\Gamma=E[X\otimes X]$, where $X\otimes X$ indicates the (random) operator such that $(X\otimes X)f=\langle X,f\rangle X$, $\forall f\in L^2(\calT)$.  
Assuming $E\|X\|^2<\infty$, for any orthonormal system $\{\varphi_j\}$ in $L^2(\calT)$, we have
\bse
\sum_j\langle\Gamma\varphi_j,\varphi_j\rangle=E[\sum_j\langle X,\varphi_j\rangle^2]=E\|X\|^2,
\ese
and thus $\Gamma$ is a trace-class operator. In particular, it is a compact operator and thus we have the spectral decomposition
\bse
\Gamma=\sum_{j=1}^\infty s_j(e_j\otimes e_j),
\ese
for some (eigenvalues) $s_1\ge s_2\ge\cdots\ge 0$ and (eigenfunctions) $e_j\in L^2(\calT)$ is an orthonormal system in $L^2(\calT)$.
For simplicity of exposition, we assume all eigenvalues are positive and $e_j, j=1,2,\ldots$ forms an orthonormal basis (i.e. $\overline{\rm span}\{e_j,j=1,2,\ldots\}=L^2(\calT)$). 

\begin{rmk}\label{rmk:pd}
We briefly discuss what happens if $\overline{\rm span}\{e_j,j=1,2,\ldots\}\subsetneq L^2(\calT)$. In this case, any $f\in L^2(\calT)$ can be written as $f_K+f_I$ with $f_K\in Ker(\Gamma)$ and $f_I\in Ker(\Gamma)^\bot=\overline{Im(\Gamma)}$. Then we have $E[\langle X,f\rangle^2]=\langle \Gamma f,f\rangle=\langle\Gamma f_I,f_I\rangle$. This implies $X\in \overline{Im(\Gamma)}$ almost surely and the component $f_K$ cannot be identified due to $\langle X,f_K\rangle=0$, a.s. That is, $f_K$ will be annihilated after taking an inner product with $X$. Thus, we need to consider error bounds for $P_I(\wh f-f_0)$, for example, where $P_I$ denotes the projection onto $\overline{Im(\Gamma)}$, instead of for $\wh f-f_0$. Alternatively, we can also simply assume the true function $f_0\in \overline{Im(\Gamma)}$ since its component in $Ker(\Gamma)$ cannot be identified anyway.
\end{rmk}
\begin{rmk}
We adopted the view that $X$ is a Hilbert-space-valued random element. Alternatively, one can adopt a stochastic process view that $X=\{X(t), t\in\calT\}$ such that $X(t)$ is a random variable for any $t$, and define the covariate kernel $\Gamma(s,t)=E[X(s)X(t)]$. The two views of $X$ are closely related and under mild assumptions, $\Gamma$ defined before is just the integral operator with kernel $\Gamma(s,t)$, and the two views lead to identical theoretical results.
\end{rmk}

In the rest of the paper, $C$ denotes a generic positive constant that can assume different values even on the same line. We impose the following technical assumptions. 
\begin{itemize}
\item[(A)] $E\|X\|^4<\infty$, $\Gamma$ is a positive-definite trace-class operator. Furthermore, we assume
\be\label{eqn:4th}
(E\langle X,f\rangle^4)^{1/4}\le C(E\langle X,f\rangle^2)^{1/2},\; \forall f\in L^2(\calT).
\ee
\item[(B)] The true function $f_0=\Gamma^{r}g_0$ for some $r\ge 0$, $g_0\in L^2(\calT)$ with $\|g_0\|\le C$ for some constant $C>0$. 
\item[(C)] For some $\beta\in(0,1]$, $\calN(\lambda):={\tr}((\Gamma+\lambda)^{-1}\Gamma)\le C\lambda^{-\beta}$.
\item[(D)] 
 The noise $\epsilon_i$ in the true model  satisfies $E[\epsilon_i|X_i]=0$ and $E[|\epsilon_i|^{1+\kappa}|X_i]\le v_\kappa$ for some $\kappa>0$, $v_\kappa>0$.
\end{itemize}
For bounds with exponential tail, we also need the following assumption.
\begin{itemize}
\item[(E)] The functional predictor has the representation $X=\sum_{j=1}^\infty \xi_je_j$ where $\xi_j=\langle X,e_j\rangle$ are independent of each other and $\xi_j$ is sub-Gaussian with parameter $s_j$ in the sense that $E[e^{t\xi_j^2/s_j}]\le e^{Ct^2}$, $j=1,2,\ldots$.
\end{itemize}

\begin{rmk}\label{rmk:ass}
As discussed above, $E\|X\|^2<\infty$ already implies $\Gamma$ is a trace-class operator. We also note that ${\rm tr}(\Gamma)=\sum_j\langle \Gamma e_j,e_j\rangle=\sum_js_j$. Equation \eqref{eqn:4th} was also used in \cite{hall07,caiyuan12}. It is satisfied if $X$ is a Gaussian process, but only requires the 4-th moment. 
Assumption (B) imposes a smoothness condition on the truth. 
 Indeed, writing $g_0=\sum_j\langle g_0,e_j\rangle  e_j$, then $\Gamma^r g_0=\sum_js_j^r\langle g_0,e_j\rangle e_j$. This means that high-frequency components (typically associated with larger $j$) are suppressed as $r$ becomes large, and thus $f_0$ is ``smoother" when $r$ increases. (B) can be written equivalently as $\sum_{j=1}^\infty f_{0j}^2/s_j^{2r}<\infty$, with $f_{0j}=\langle f_0,e_j\rangle.$ When $s_j\asymp j^{-\alpha}$ for some $\alpha>1$ (as often assumed in the FLR literature), a sufficient condition for (B) is that $f_{0j}\le Cj^{-b}$ for $b>\alpha r+1/2$.
  For (C), we note that if $s_j\asymp j^{-\alpha}$ for some $\alpha>1$, we will have $\calN(\lambda):={\rm tr}(\Gamma+\lambda)^{-1}\Gamma)=\sum_j\frac{s_j}{s_j+\lambda}\le C\lambda^{-\frac{1}{\alpha}}$. Furthermore, since $\calN(\lambda)=\sum_j\frac{s_j}{s_j+\lambda}\le \sum_j\frac{s_j}{\lambda}$, (C) is always satisfied with $\beta=1$. In our main result, we actually only use the definition of $\calN(\lambda)$, and the bound $C\lambda^{-\beta}$ is used in Remark \ref{rmk:Nlambda} to simplify the expressions of the main result.
  (D) is the main assumption that differs from typical assumptions in the least squares case. Unlike the least squares case, we do not need the second-moment assumption on errors to establish some statistical rate (although the rate is slower if $\kappa<1$) and we do not need an exponential tail bound for the error distribution in order to establish the estimator's convergence rate with an exponential tail bound. $v_\kappa$ can usually be regarded as a constant, but we will still track the effect of $v_\kappa$ explicitly in our bounds. (E) is satisfied if $X$ is Gaussian. Since $e_j, j=1,2,\ldots,$ are the eigenfunctions of $\Gamma$, we always have $E[\xi_j\xi_k]=E\langle \Gamma e_j,e_k\rangle=0$ for $j\neq k$. That is, $\xi_j$ and $\xi_k$ are uncorrelated. If $X$ is Gaussian, then $\xi_j$ and $\xi_k$ are independent. Furthermore, $E[\xi_j^2]=s_j$ and thus it is natural to standardize $\xi_j^2$ by $s_j$. 
\end{rmk}
We can now state our main result. Let $\check\kappa=\min\{\kappa ,1\}$ and $\check r=\min\{r,1/2\}$. Our rate is somewhat innovatively obtained for $\|(\Gamma+\lambda)^{1/2}(\wh f-f_0)\|$, which seems to be new in the literature. It is easy to see that $\|(\Gamma+\lambda)^{1/2}(\wh f-f_0)\|\asymp \|\Gamma^{1/2}(\wh f-f_0)\|+\sqrt{\lambda}\|\wh f-f_0\|$, and thus we can actually obtain the rate for both the prediction error $\|\Gamma^{1/2}(\wh f-f_0)\|$ and for $L^2$ error $\|\wh f-f_0\|$ in a \emph{unified} way. We also note that $\|\Gamma^{1/2}f\|^2=E[\langle X,f\rangle^2]$, which is the reason we call it prediction error.
\begin{thm}\label{thm:mainflr}
Define $\omega_n=\sqrt{v_{\check\kappa}\gamma^{1-\check\kappa}\frac{\calN(\lambda)}{n}}+\lambda^{\check{r}+1/2}+v_{\check\kappa}\gamma^{-\check\kappa}$. Under Assumptions (A)-(D), and that $\gamma\ge 2(4v_{\check\kappa})^{\frac{1}{1+\check\kappa}}$, $\gamma\ge C\omega_n$ for a sufficiently large $C$, $\frac{\calN(\lambda)}{n\lambda}=o(1)$,
we have
\bse
\|(\Gamma+\lambda)^{1/2}(\wh f-f_0)\|=O_p(\omega_n).
\ese
If in addition (E) holds, and $C\omega_n\log(\frac{n}{\delta})\le \gamma \le Cn^{\frac{1}{1+\kappa}}$ for a sufficiently large $C$, $\frac{\calN(\lambda)}{n\lambda}\log^2(\frac{n}{\delta})=o(1)$,  we have with probability $1-\delta$,
\bse
\|(\Gamma+\lambda)^{1/2}(\wh f-f_0)\|\le C\omega_n\log(\frac{n}{\delta}).
\ese
\begin{rmk} The three terms in $\omega_n$ correspond to the stochastic error, the bias  due to regularization, and the bias due to using Huber's loss, respectively. 
The rate improves as $\kappa$ and $r$ increases, but does not further improve beyond $\kappa=1$ and $r=1/2$. This ``saturation effect" for $\kappa$ is the same as in parameter models using the adaptive Huber's loss. The ``saturation effect" for $r=1/2$ is similar to that appeared in kernel ridge regression \citep{caponnetto07}. 
\end{rmk}
\begin{rmk}\label{rmk:Nlambda} To see the rate more explicitly, we assume $v_\kappa$ is a constant and set $\lambda$, $\gamma$ such that the three terms in $\omega_n$ are of the same order. This leads to $\lambda\asymp n^{-\frac{2\kappa}{(2r+1)(1+\kappa)+2\beta\kappa}}$, $\gamma\asymp n^{\frac{2r+1}{(2r+1)(1+\kappa)+2\beta\kappa}}$, and $\omega_n\asymp n^{-\frac{(2r+1)\kappa}{(2r+1)(1+\kappa)+2\beta\kappa}}$, for $\kappa\in(0,1]$, $r\in[0,1/2]$. The rate in $L^2$ is $\|\wh f-f_0\|=O_p(n^{-\frac{2r\kappa}{(2r+1)(1+\kappa)+2\beta\kappa}})$ (assumptions regarding $\gamma,\lambda$ stated in Theorem \ref{thm:mainflr} are automatically satisfied for these specific choices of $\gamma$, $\lambda$).

When $\kappa=1$, the rate obtained is the same as that for least squares regression (for bounds with exponential tail, there is an extra factor of $\log n$) obtained in \cite{Zhang2020}, which extends the results of \cite{caiyuan12} to the case $r>0$. 
\end{rmk}

\end{thm}
\noindent\textbf{Proof of Theorem \ref{thm:mainflr}.}
We prove this by contraction. To prove the rate in probability, assume that there exists some constant $\delta>0$ such that 
$P(\|(\Gamma+\lambda)^{1/2}(\wh f-f_0)\|>L\omega_n)>\delta$ for all $L$ (however large it is).
By choosing $v=\frac{L\omega_n}{\max\{\|(\Gamma+\lambda)^{1/2}(\wh f-f_0)\|,L\omega_n\}}\in(0,1]$ and defining $f=v\wh f+(1-v)f_0$, this $f$ satisfies 
\be\label{eqn:contra1}
\|(\Gamma+\lambda)^{1/2}( f-f_0)\|\le L\omega_n \mbox{ and } P(\|(\Gamma+\lambda)^{1/2}( f-f_0)\|=L\omega_n)>\delta.
\ee

To prove the exponential tail bound, we assume that for any $L>0$, there exists some $\delta\in (0,1)$ such that 
\bse
P(\|(\Gamma+\lambda)^{1/2}( f-f_0)\|>L\omega_n\log(n/\delta))>\delta,
\ese
and with $v=\frac{L\omega_n\log(n/\delta)}{\max\{\|(\Gamma+\lambda)^{1/2}(\wh f-f_0)\|,L\omega_n\log(n/\delta)\}}\in(0,1]$, $f=v\wh f+(1-v)f_0$, we have that
\be\label{eqn:contra2}
P(\|(\Gamma+\lambda)^{1/2}( f-f_0)\|=L\omega_n\log(n/\delta))>\delta.
\ee
In the following, $f$ without other qualifications always indicates this particular $f$ (slightly different definitions of $f$ for bound in probability and bound with exponential tail).

\noindent\textbf{Step 1}. We show that  $
\langle \nabla H_\gamma(  f)-\nabla H_\gamma(f_0),  f-f_0\rangle \le  \langle \nabla H_\gamma( \wh f)-\nabla H_\gamma(f_0),    f-f_0\rangle.$

Convexity of $H_\gamma$ implies that for any $f_1,f_2\in L^2(\calT)$,
\be\label{eqn:conv0}
\langle \nabla H_\gamma(f_1)-\nabla H_\gamma(f_2),f_1-f_2\rangle\ge 0.
\ee

 Thus we have 
\be\label{eqn:fff}
\langle \nabla H_\gamma(\wh f)-\nabla H_\gamma(f),\wh f-f_0\rangle=\frac{1}{1-v}\langle \nabla H_\gamma(\wh f)-\nabla H_\gamma(f),\wh f-f\rangle\ge 0,
\ee
where the equality used that $\wh f-f=(1-v)(\wh f-f_0)$ by the definition of $v$ and $f$, and the inequality is due to \eqref{eqn:conv0}.
Thus,
\bse
\langle \nabla H_\gamma(  f)-\nabla H_\gamma(f_0),  f-f_0\rangle&=&v\langle \nabla H_\gamma(  f)-\nabla H_\gamma(f_0),  \wh f-f_0\rangle\le v\langle \nabla H_\gamma( \wh f)-\nabla H_\gamma(f_0),  \wh f-f_0\rangle\\
&=&\langle \nabla H_\gamma( \wh f)-\nabla H_\gamma(f_0),  f-f_0\rangle,
\ese
where the first step used $f-f_0=v(\wh f-f_0)$, the second step used \eqref{eqn:fff}, and the last step again used $f-f_0=v(\wh f-f_0)$.

\noindent\textbf{Step 2.} We give a lower bound for $P_n \langle \nabla H_\gamma(  f)-\nabla H_\gamma(f_0),  f-f_0\rangle$.

Define $\calE_i=\{|\epsilon_i|\le \frac \gamma 2, |\langle X_i, f-f_0\rangle|\le \frac \gamma 2\}$, $\calE=\{|\epsilon|\le \frac \gamma 2, |\langle X, f-f_0\rangle|\le \frac \gamma 2\}$.
We have 
\bse
&&P_n\langle \nabla H_\gamma(  f)-\nabla H_\gamma(f_0),  f-f_0\rangle\\
&\ge& P_n\langle \nabla H_\gamma(  f)-\nabla H_\gamma(f_0),  f-f_0\rangle I_\calE\\
&=&E[\langle \nabla H_\gamma(  f)-\nabla H_\gamma(f_0),  f-f_0\rangle I_\calE]+(P_n-P)\langle \nabla H_\gamma(  f)-\nabla H_\gamma(f_0),  f-f_0\rangle I_\calE.
\ese
 On $\calE$, we have  $ \langle\nabla H_\gamma( f)- \nabla H_\gamma( f_0),f-f_0\rangle=\langle X,f-f_0\rangle^2$.  
Thus, we get
\bse
&& E \langle \nabla H_\gamma( f)-\nabla H_\gamma(f_0),f-f_0\rangle I_{\calE}\\
&=& E[\langle X,f-f_0\rangle^2I_{\calE}]\\
&=&E[ \langle X,f-f_0\rangle^2] -E[ \langle X,f-f_0\rangle^2 I_{\calE^c}]\\
&\ge &E[ \langle X,f-f_0\rangle^2]-E[ \langle X,f-f_0\rangle^2]v_\kappa(\frac{2}{\gamma})^{1+\kappa}-C(E[ \langle X,f-f_0\rangle^4])^{1/2}\frac{2\|\Gamma^{1/2}(f-f_0)\|}{\gamma}\\
&\ge &\frac 1 2 E[ \langle X,f-f_0\rangle^2],
\ese
where the penultimate step used Markov's inequality $P(|\ep|>\frac{\gamma}{2}|x)\le \nu_\kappa(2/\gamma)^{1+\kappa}$ and $P(\langle X,f-f_0\rangle>\frac{\gamma}{2})\le \frac{4E[\langle X,f-f_0\rangle^2]}{\gamma^2}$, and the final inequality is due to our assumption on $\gamma$ in the statement of the theorem ($\gamma$ is sufficiently large).

Then we consider $(P_n-P)[\langle X,f-f_0\rangle^2I_{\calE}]$.
We have that, denoting $g=(\Gamma+\lambda)^{1/2}(f-f_0)$
\bse
&&(P_n-P)[\langle X,f-f_0\rangle^2I_{\calE}]\\
&=&(P_n-P)\langle (\Gamma+\lambda)^{-1/2}XI_{\calE}, g \rangle^2I_{\calE}\\
&=&\langle(\Gamma+\lambda)^{-1/2}(\wt\Gamma_n-\wt\Gamma)(\Gamma+\lambda)^{-1/2}g,g\rangle\\
&\le &\|(\Gamma+\lambda)^{-1/2}(\wt\Gamma_n-\wt\Gamma)(\Gamma+\lambda)^{-1/2}\|_{HS}\|g\|^2\\
&\le &\|(\Gamma+\lambda)^{-1/2}(\wt\Gamma_n-\wt\Gamma)\|_{HS}\cdot\frac{1}{\sqrt{\lambda}}\|(\Gamma+\lambda)^{1/2}(f-f_0)\|^2,
\ese
where $\wt\Gamma_n=\frac{1}{n}\sum_i\wt X_i\otimes \wt X_i$, $\wt X_i=X_iI_{\calE_i}$, and $\wt\Gamma=E[\wt X_i\otimes \wt X_i]$.
We now use some technical arguments to remove the indicator $I_{\calE}$. More specifically, for rates in probability, we use
\bse
&&E\|(\Gamma+\lambda)^{-1/2}(\wt\Gamma_n-\wt\Gamma) \|_{HS}\\
&\le &E\|\frac{1}{n}\sum_i\sigma_i (\Gamma+\lambda)^{-1/2}(X_i\otimes X_i)I_{\calE_i} \|_{HS}\\
&\le &E\|\frac{1}{n}\sum_i\sigma_i (\Gamma+\lambda)^{-1/2}(X_i\otimes X_i) \|_{HS},
\ese
where the first step used symmetrization (Theorem 2.1 of \cite{koltchinskiibook2011}), with $\sigma_i\in\{-1,1\}$ being i.i.d. Rademacher variables independent of all other random variables, and the second step used the contraction inequality to remove the indicator (Theorem 4.4 of \cite{ledoux91}).
For convergence with exponential tail bound, we use that, for any $a>4E\|(\Gamma+\lambda)^{-1/2}(\wt\Gamma_n-\wt\Gamma) \|_{HS}$, 
\bse
&&P(\|(\Gamma+\lambda)^{-1/2}(\wt\Gamma_n-\wt\Gamma)\|_{HS}>a)\\
&\le&\frac{2P(\|\frac{1}{n}\sum_i(\Gamma+\lambda)^{-1/2}\sigma_iI_{\calE_i}(X_i\otimes X_i)\|_{HS}>a/4)}{1-\frac{2E\|(\Gamma+\lambda)^{-1/2}(\wt\Gamma_n-\wt\Gamma) \|_{HS}}{a}}\\
&\le& 4P(\|\frac{1}{n}\sum_i(\Gamma+\lambda)^{-1/2}\sigma_iI_{\calE_i}(X_i\otimes X_i)\|_{HS}>a/4)\\
&\le &8P(\|\frac{1}{n}\sum_i(\Gamma+\lambda)^{-1/2}\sigma_i (X_i\otimes X_i)\|_{HS}>a/4),
\ese
where the first step used symmetrization for probabilities (Lemma 2.3.7 of \cite{vaartwellner96}),  the second step used $a>4E\|(\Gamma+\lambda)^{-1/2}(\wt\Gamma_n-\wt\Gamma) \|_{HS}$, and the last step used another version of contraction inequality for probabilities (Theorem 4.4 of \cite{ledoux91}) to remove the indicator.

Thus, the bounds (in probability and with exponential tail) for $\|(\Gamma+\lambda)^{-1/2}(\wt\Gamma_n-\wt\Gamma)\|_{HS}$ are the same as in Lemma \ref{prp:ttn} for $\|(\Gamma+\lambda)^{-1/2}(\Gamma_n- \Gamma)\|_{HS}$.

\noindent\textbf{Step 3.} We give an upper bound for 
$ P_n\langle \nabla H_\gamma( \wh f)-\nabla H_\gamma(f_0), f-f_0\rangle$.

Since $\wh f$ satisfies the first-order optimality condition of the objective function, we have 
\bse
&&P_n\langle \nabla H_\gamma( \wh f), f-f_0\rangle=-\lambda\langle \wh f,f-f_0\rangle=-\frac{\lambda}{v}\|f-f_0\|^2-\lambda\langle f_0,f-f_0\rangle\\
&\le & -\frac{\lambda}{v}\|f-f_0\|^2 +C\lambda^{r+\frac{1}{2}}(\|\Gamma^{1/2}(f-f_0)\|+\sqrt{\lambda}\|f-f_0\|),
\ese
where we used again $f-f_0=v(\wh f-f_0)$ and that 
\bse
&&\lambda|\langle f_0,f-f_0\rangle|\\
&=&\lambda|\langle \Gamma^r g_0,f-f_0\rangle|\\
&=&\lambda|\langle  g_0,\Gamma^r(f-f_0)\rangle|\\
&\le &C\lambda^{r+\frac{1}{2}}\|\lambda^{\frac{1}{2}-r}\Gamma^r(f-f_0)\|\\
&= &C\lambda^{r+\frac{1}{2}}\langle f-f_0, \lambda^{1-2r}\Gamma^{2r}(f-f_0)\rangle^{1/2}\\
&\le & C\lambda^{r+\frac{1}{2}}\langle f-f_0,((1-2r)\lambda+2r\Gamma)(f-f_0)\rangle^{1/2}\\
&=&C\lambda^{r+\frac{1}{2}}\left((1-2r)\lambda\|f-f_0\|^2+2r\|\Gamma^{1/2}(f-f_0)\|^2\right)^{1/2},
\ese
where the second inequality used Young's inequality for operators: $\lambda^{1-2r}\Gamma^{2r}\le (1-2r)\lambda+2r\Gamma$.

Furthermore, we have
\be\label{eqn:error2}
&&P_n\langle \nabla H_\gamma(f_0),f-f_0\rangle\nonumber\\
&=&P_n\langle \nabla E[H_\gamma(f_0)|X],f-f_0\rangle+P_n\langle \nabla H_\gamma(f_0)-E[\nabla H_\gamma(f_0)|X],f-f_0\rangle.
\ee
For the first term on the right-hand side of \eqref{eqn:error2}, we have
\bse
&&|P_n\langle  E[\nabla H_\gamma(f_0)|X],f-f_0\rangle|\\
&=&|P_n\langle  E[H_\gamma'(\epsilon)|X]X,f-f_0\rangle|\\
&\le & v_\kappa\gamma^{-\kappa} P_n\langle X,f-f_0\rangle\\
&\le &v_\kappa\gamma^{-\kappa} (P_n\langle X,f-f_0\rangle^2)^{1/2}\\
&=& v_\kappa\gamma^{-\kappa} (\langle\Gamma_n(f-f_0),f-f_0\rangle)^{1/2},
\ese
where the first inequality above used Lemma \ref{lem:moment}. Thus 
\bse
&&(P_n\langle  E[\nabla H_\gamma(f_0)|X],f-f_0\rangle)^2\\
&\le &v_\kappa^2\gamma^{-2\kappa}\|(\Gamma+\lambda)^{-1/2}\Gamma_n(\Gamma+\lambda)^{-1/2}\|_{op}\cdot\|(\Gamma+\lambda)^{1/2}(f-f_0)\|^2.
\ese 
 For the 2nd term on the right-hand side of \eqref{eqn:error2}, we use
\be\label{eqn:eta}
&&|P_n\langle \nabla H_\gamma(f_0)-E[\nabla H_\gamma(f_0)|X],f-f_0\rangle|\nonumber\\
&\le &\|P_n( (\Gamma+\lambda)^{-1/2}X\underbrace{(H_\gamma'(\epsilon)-E[H_\gamma'(\epsilon)|X])}_{\eta})\|\cdot \|(\Gamma+\lambda)^{1/2}(f-f_0)\|,
\ee
and then we can apply the bound from Lemma \ref{lem:ep}.

\noindent\textbf{Step 4.} Combining Steps 1-3,
we get
\be\label{eqn:get1}
&&\frac{1}{2}\|\Gamma^{1/2}(f-f_0)\|^2+\frac{\lambda}{v}\|f-f_0\|^2\nonumber\\
&\le &\frac{1}{\sqrt{\lambda}}\|(\Gamma+\lambda)^{-1/2}(\wt\Gamma_n-\wt\Gamma)\|_{HS}\cdot\|\Gamma^{1/2}(f-f_0)\|^2\nonumber\\
&&+C\lambda^{r+1/2}\|(\Gamma+\lambda)^{1/2}(f-f_0)\| \nonumber\\
&&+v_\kappa\gamma^{-\kappa}\|(\Gamma+\lambda)^{-1/2}\Gamma_n(\Gamma+\lambda)^{-1/2}\|_{op}^{1/2}\|(\Gamma+\lambda)^{1/2}(f-f_0)\|\nonumber\\
&&+\|P_n( (\Gamma+\lambda)^{-1/2}X {\eta})\|\cdot \|(\Gamma+\lambda)^{1/2}(f-f_0)\|\nonumber\\
&=&\frac{1}{\sqrt{\lambda}}\|(\Gamma+\lambda)^{-1/2}(\wt\Gamma_n-\wt\Gamma)\|_{HS}\cdot\|\Gamma^{1/2}(f-f_0)\|^2\nonumber\\
&&+C\left(\lambda^{r+1/2}+\|P_n( (\Gamma+\lambda)^{-1/2}X {\eta})\|+v_\kappa\gamma^{-\kappa}\|(\Gamma+\lambda)^{-1/2}\Gamma_n(\Gamma+\lambda)^{-1/2}\|_{op}^{1/2}\right)\nonumber\\
&&\hspace{.1in}\cdot\|(\Gamma+\lambda)^{1/2}(f-f_0)\|.
\ee
Bound for the term $\|(\Gamma+\lambda)^{-1/2}(\wt\Gamma_n-\wt\Gamma)\|_{HS}$ has been discussed in Step 2, and for $\|(\Gamma+\lambda)^{-1/2}\Gamma_n(\Gamma+\lambda)^{-1/2}\|_{op}$, we can use
\be\label{eqn:get2}
&&\|(\Gamma+\lambda)^{-1/2}\Gamma_n(\Gamma+\lambda)^{-1/2}\|_{op}\nonumber\\
&\le &\|(\Gamma+\lambda)^{-1/2}\Gamma(\Gamma+\lambda)^{-1/2}\|_{op}+\|(\Gamma+\lambda)^{-1/2}(\Gamma_n-\Gamma)(\Gamma+\lambda)^{-1/2}\|_{op}\nonumber\\
&\le & 1+\|(\Gamma+\lambda)^{-1/2}(\Gamma_n-\Gamma)\|_{HS}\lambda^{-1/2},
\ee
and then apply Lemma \ref{prp:ttn}.

More specifically, to get the bound in probability, we get from \eqref{eqn:get1} and \eqref{eqn:get2} that 
\bse
&&\left(1-O_p\left(\sqrt{\frac{\calN(\lambda)}{n\lambda}}\right)\right)\|(\Gamma+\lambda)^{1/2}(f-f_0)\|^2\\
&=&O_p\left(\sqrt{v_\kappa\gamma^{1-\kappa}\frac{\calN(\lambda)}{n}}+\lambda^{r+1/2}+v_\kappa\gamma^{-\kappa}\sqrt{1+\sqrt{\frac{\calN(\lambda)}{n\lambda}}}\right)\|(\Gamma+\lambda)^{1/2}(f-f_0)\|,
\ese
which implies
\bse
\|(\Gamma+\lambda)^{1/2}(f-f_0)\|=O_p(\omega_n).
\ese
This implies that $P(\|(\Gamma+\lambda)^{1/2}(f-f_0)\|\ge L\omega_n)<\delta$ for $L$ sufficiently large, which contradicts \eqref{eqn:contra1}.

To get the bound with exponential tail probability, we also have from \eqref{eqn:get1} and \eqref{eqn:get2} that with probability $1-\delta$ (note $\gamma=O(n^{\frac{1}{1+\kappa}})$ ensures $\sqrt{v_\kappa\gamma^{1-\kappa}\frac{\calN(\lambda)}{n}}\log(\frac{n}{\delta})$ is the dominating term in \eqref{eqn:ep22} of Lemma \ref{lem:ep}),
\bse
&&\left(1- C\sqrt{\frac{\calN(\lambda)}{n\lambda}}\log(\frac{n}{\delta})\right) \|(\Gamma+\lambda)^{1/2}(f-f_0)\|\\
&\le&C\left(\sqrt{v_\kappa\gamma^{1-\kappa}\frac{\calN(\lambda)}{n}}\log(\frac{n}{\delta})+\lambda^{r+1/2}+v_\kappa\gamma^{-\kappa}\sqrt{1+\sqrt{\frac{\calN(\lambda)}{n\lambda}}\log(\frac{n}{\delta})}\right), 
\ese
which contradicts \eqref{eqn:contra2}.
\hfill $\Box$

\begin{lem}\label{prp:ttn}
Under Assumption (A), we have
\be
E\|(\Gamma+\lambda I)^{-1/2}(\Gamma-\Gamma_n)\|_{HS}^2&=&O(\frac{\calN(\lambda)}{n}).\label{eqn:ttn1}
\ee
Furthermore, if Assumption (E) also holds, then for any $0<\delta<1$, with probability at least $1-\delta$,
\begin{equation}\label{eqn:ttn3}
\|(\Gamma+\lambda I)^{-1/2}(\Gamma-\Gamma_n)\|_{HS}\le C\left(\frac{\sqrt{M_1(\delta)M_2(\delta)}}{n}+\sqrt{\frac{ {\calN(\lambda)}}{n}}\right)\log (\frac{n}{\delta})+C\sqrt{\calN(\lambda)\frac{\delta}{n}},
\end{equation}
where $M_1(\delta)=\sum_js_j^2+\log(\frac{n}{\delta})$ and $M_2(\delta)=\calN(\lambda)+\log(\frac{n}{\delta})$.

The above bounds still hold when $\Gamma_n-\Gamma$ is replaced by $\Gamma_{n,\sigma}:=\frac{1}{n}\sum_i\sigma_i(X_i\otimes X_i)$, where $\sigma_i\in\{-1,1\}$ are i.i.d. Rademacher variables independent of all other random variables. 
\end{lem} 
\textbf{Proof of Lemma \ref{prp:ttn}}.
For a random predictor $X=\sum_{j=1}^\infty\xi_je_j$,
\be\label{eqn:infbound}
&&\|(\Gamma+\lambda)^{-1/2}(X\otimes X)\|_{HS}^2\nonumber\\
&=&\sum_{j,k}\langle(\Gamma+\lambda)^{-1/2}(X\otimes X)e_j,e_k\rangle^2\nonumber\\
&=&\sum_{j,k}\xi_j^2\langle(\Gamma+\lambda)^{-1/2}X,e_k\rangle^2\nonumber\\
&=&\sum_{j,k}\xi_j^2\langle X,(\Gamma+\lambda)^{-1/2}e_k\rangle^2\nonumber\\
&=&\sum_{j,k}\xi_j^2\langle X,e_k/\sqrt{s_k+\lambda}\rangle^2\nonumber\\
&=&\sum_{j,k}\frac{\xi_j^2\xi_k^2}{(s_k+\lambda)}.
\ee 
This implies
\bse
E\|(\Gamma+\lambda)^{-1/2}(X\otimes X)\|_{HS}^2\le C\sum_{j,k}\frac{s_js_k}{(s_k+\lambda)}\le C\sum_{k}\frac{s_k}{(s_k+\lambda)}=C\calN(\lambda).
\ese

Furthermore, under Assumption (E), since $\xi_j/\sqrt{s_j}$ is sub-Gaussian, we know that $\xi_j^2/s_j$ is sub-exponential (Lemma 2.7.7 of \cite{vershynin18book}). Then by Proposition 2.7.1 (e) of \cite{vershynin18book},  there exist some constants $c_0,c_1>0$ such that
\bse
E[\exp\{t(\xi_k^2/s_k)\}]\le \exp\{c_0t^2\}, |t|\le c_1.
\ese 
Thus we have
\be\label{eqn:c1}
&& E[\exp\{t(\sum_k\frac{\xi_k^2}{s_k+\lambda})\}]\nonumber\\
&= & E[\exp\{\sum_k \frac{ts_k}{s_k+\lambda}\frac{\xi_k^2}{s_k}\}]\nonumber\\
&\le & \exp\{c_0t^2\sum_k\frac{s_k^2}{(s_k+\lambda)^2}\}, |t|\le c_1,
\ee
and thus
\bse
&&P(\sum_k\frac{\xi_k^2}{s_k+\lambda}>a)\\
&\le & E[\exp\{t(\sum_k\frac{\xi_k^2}{s_k+\lambda})\}]/e^{at}\\
&\le & \exp\{c_0t^2\calN(\lambda)-at\}, |t|\le c_1.
\ese
In particular, plugging in $t=c_1$ above we get
\bse
&&P(\sum_k\frac{\xi_k^2}{s_k+\lambda}>a)\\
&\le & \exp\{c_0c_1^2\calN(\lambda)-c_1a\}.
\ese
Letting $a=c_0c_1\calN(\lambda)+c_1^{-1}\log(\delta^{-1})$, we have
\be\label{eqn:delta1}
P(\sum_k\frac{\xi_k^2}{s_k+\lambda}>c_0c_1\calN(\lambda)+c_1^{-1}\log(\delta^{-1}) )&\le & \delta.
\ee
Similarly, using 
\bse
&&E[\exp\{t\sum_j\xi_j^2\}]\\
&\le & \exp\{c_0t^2\sum_js_j^2\},\; |t|\le c_2:=c_1/s_1,
\ese
we have for all $a>0$
\bse
P(\sum_j \xi_j^2>a)\le \exp\{c_0t^2\sum_js_j^2-at\},\; |t|\le c_2.
\ese
Setting $t=c_2$ and $a=c_0c_2\sum_j s_j^2+c_2^{-1}\log(\delta^{-1})$,
\be\label{eqn:delta2}
P(\sum_j \xi_j^2>c_0c_2\sum_j s_j^2+c_2^{-1}\log(\delta^{-1}))\le \delta.
\ee

Combining \eqref{eqn:infbound},  \eqref{eqn:delta1} and \eqref{eqn:delta2}, we get
\be\label{eqn:delta3}
&&P\Big(\|(\Gamma+\lambda I)^{-1/2}(X\otimes X)\|_{HS}^2>(c_0c_2\sum_j s_j^2+c_2^{-1}\log(2\delta^{-1}))\nonumber\\
&&\hspace{2.2in}( c_0c_1\calN(\lambda)+c_1^{-1}\log(2\delta^{-1})) \Big)\le \delta.
\ee
Let $M_1(\delta):=c_0c_2\sum_j s_j^2+c_2^{-1}\log(2\delta^{-1})$, $M_2(\delta):=c_0c_1\calN(\lambda)+c_1^{-1}\log(2\delta^{-1})$, 
$$\Omega(\delta)=\{\sum_j\xi_j^2\le M_1(\delta), \sum_j\xi_j^2/(s_j+\lambda)\le  M_2(\delta)\},$$ and define $\wt X=XI\{\Omega(\delta)\}$.
Then \eqref{eqn:infbound} implies
$$\|(\Gamma+\lambda I)^{-1/2}(\wt X\otimes \wt X)\|_{HS}^2\le M_1(\delta)M_2(\delta).
$$
Furthermore, $E\|(\Gamma+\lambda I)^{-1/2}(\wt X\otimes \wt X)\|_{HS}^2\le E\|(\Gamma+\lambda I)^{-1/2}(X\otimes X)\|_{HS}^2\le C\calN(\lambda)$, and by Lemma 24 of \cite{Lin2020} (Bernstein's inequality for random elements in a Hilbert space), with probability at least $1-\delta$,
\be\label{eqn:comb1}\|(\Gamma+\lambda I)^{-1/2}(\wt\Gamma-\wt\Gamma_n)\|_{HS}\le C\left(\frac{\sqrt{M_1(\delta)M_2(\delta)}}{n}+\sqrt{\frac{{\calN(\lambda)}}{n}}\right)\log (2\delta^{-1}),
\ee
where $\wt\Gamma_n=(1/n)\sum_i(\wt X_i\otimes \wt X_i)$ and $\wt\Gamma=E[\wt\Gamma_n]$.

We have obviously
\be\label{eqn:comb2}
P(\|(\Gamma+\lambda I)^{-1/2}(\Gamma_n-\wt\Gamma_n)\|_{HS}\neq 0)\le nP(\Omega^{c}(\delta))\le n\delta.
\ee
Using $\Gamma-\wt \Gamma=E[I(\Omega^{c}(\delta))(X\otimes X)]$, we have
\be\label{eqn:comb3}
&&\|(\Gamma+\lambda I)^{-1/2}(\Gamma-\wt\Gamma)\|_{HS}^2\nonumber\\
&=& \sum_{j,k} \langle E[I(\Omega^{c}(\delta))(X\otimes X)]e_j,\frac{e_k}{\sqrt{s_k+\lambda}}\rangle^2\nonumber\\
&=& \sum_{j,k}\left(E\left[\frac{\xi_j\xi_k}{\sqrt{s_k+\lambda}}I(\Omega^{c}(\delta))\right]\right)^2\nonumber\\
&\le&\sum_{j,k}\frac{s_js_k}{(s_k+\lambda)}P(\Omega^{c}(\delta))\nonumber\\
&\le & C\calN(\lambda)P(\Omega^{c}(\delta))\nonumber\\
&\le&  C\calN(\lambda)\delta,
\ee
where the first inequality above used the H\"older's inequality. Combining \eqref{eqn:comb1}-\eqref{eqn:comb3}, and replacing $\delta$ with $\delta/(n+1)$ (to make sure the bound holds with probability $1-\delta$), we proved the bound \eqref{eqn:ttn3}. 

Finally, for $\Gamma_{n,\sigma}$, it is easy to see that the proofs still go through without change, noting that $\|(\Gamma+\lambda)^{-1/2}\sigma(X\otimes X)\|_{HS}=\|(\Gamma+\lambda)^{-1/2}(X\otimes X)\|_{HS}$.
\hfill $\Box$
\begin{lem}\label{lem:ep}
Under Assumptions (A) and (D), with $\eta_i=H_\gamma'(\epsilon_i)-E[H_\gamma'(\epsilon_i)|X_i]$ as defined in \eqref{eqn:eta}, we have
\bse
E\|(\Gamma+\lambda I)^{-1/2}\sum_iX_i\eta_i/n\|^2\le Cv_\kappa \gamma^{1-\kappa}\frac{\calN(\lambda)}{n}.
\ese
Furthermore, if in addition Assumption (E) holds, then
for any $0<\delta<1$, with probability at least $1-\delta$,
\be\label{eqn:ep22}
\|(\Gamma+\lambda I)^{-1/2}\sum_iX_i\eta_i/n\|\le C\left(\frac{\sqrt{M(\delta)}}{n}+\sqrt{\frac{v_\kappa\gamma^{1-\kappa}\calN(\lambda)}{n}}\right)\log(\frac{n}{\delta})+C\sqrt{v_\kappa\gamma^{1-\kappa}\calN(\lambda)\frac{\delta}{n}},
\ee
where $M(\delta)=\gamma^2(\calN(\lambda)+\log(\frac{n}{\delta}))$.
\end{lem}
\noindent\textbf{Proof of Lemma \ref{lem:ep}.} The proof is similar to that for Lemma \ref{prp:ttn}. We have
\be\label{eqn:Xeta}
\|(\Gamma+\lambda I)^{-1/2}X\eta\|^2&=&\eta^2\sum_k\xi_k^2/(s_k+\lambda),
\ee
and thus, using Lemma \ref{lem:moment},
\bse
E\|(\Gamma+\lambda I)^{-1/2}X\eta\|^2\le Cv_\kappa \gamma^{1-\kappa}\calN(\lambda).
\ese
\eqref{eqn:Xeta} also implies, using $|\eta|\le\gamma$ and \eqref{eqn:c1},
\bse
&&E[\exp\{t\|(\Gamma+\lambda I)^{-1/2}X\eta\|^2\}]\\
&=&E[E[\exp \{t\eta^2\}|X]\exp\{\sum_k\xi_k^2/(s_k+\lambda)\}]\\
&\le &E[\exp\{\gamma^2 t\sum_k\xi_k^2/(s_k+\lambda)\}]\\
&\le &\exp\{c_0\gamma^4t^2\calN(\lambda)\}, |t|\le \frac{c_1}{\gamma^2}.
\ese
Thus by Markov's inequality
$$
P(\|(\Gamma+\lambda I)^{-1/2}X\eta\|^2>a)\le \exp\{c_0\gamma^4t^2\calN(\lambda)-at\}, |t|\le \frac{c_1}{\gamma^2}.
$$
With $t=\frac{c_1}{\gamma^2}$ and $a=c_0\gamma^4t\calN(\lambda)+\frac{\gamma^2}{c_1}\log(\delta^{-1})$, we get
$$
P(\|(\Gamma+\lambda I)^{-1/2}X\eta\|^2>M(\delta))\le \delta.
$$
Defining $\Omega(\delta)=\{\|(\Gamma+\lambda I)^{-1/2}X\eta\|^2\le M(\delta)\}$ and $\wt{X\eta}=X\eta I\{\Omega(\delta)\}$, by Bernstein's inequality, we have with probability at least $1-\delta$, 
$$\|(\Gamma+\lambda I)^{-1/2}(E[\wt{X\eta}]-\sum_i\wt{X_i\eta_i}/n)\|^2\|\le C\left(\frac{\sqrt{M(\delta)}}{n}+\sqrt{\frac{v_\kappa \gamma^{1-\kappa}\calN(\lambda)}{n}}\right)\log(\delta^{-1}).
$$
Furthermore, we have
$$
P(\|(\Gamma+\lambda I)^{-1/2}\sum_i(X_i\eta_i-\wt{X_i\eta_i})/n\|\neq 0)\le n\delta,
$$
and
\bse
&&\|(\Gamma+\lambda I)^{-1/2}E[\wt{X\eta}]\|^2\\
&=&\sum_k\left(E\left[\frac{\xi_k\eta}{\sqrt{s_k+\lambda}}I\{\Omega^c(\delta)\}\right]\right)^2\\
&\le& Cv_\kappa \gamma^{1-\kappa}\sum_k\frac{s_k}{s_k+\lambda}P(\Omega^c(\delta))\\
&\le & Cv_\kappa \gamma^{1-\kappa}\calN(\lambda)\delta,
\ese
which proved the lemma.
\hfill $\Box$

\begin{lem}\label{lem:moment} Under Assumption (D), we have
\bse
&&|E[H_\gamma'(\ep_i)|X_i]|\le v_\kappa \gamma^{-\kappa}\\
&&E[(H_\gamma'(\ep_i))^2|X_i]\le 2v_\kappa \gamma^{1-\kappa}.
\ese
\end{lem}
\noindent\textbf{Proof of Lemma \ref{lem:moment}.}  
Since
\bse
H_\gamma'(\ep_i)&=&\epsilon_iI\{|\epsilon_i|\le\gamma\}+\gamma I\{\epsilon_i>\gamma\}-\gamma I\{\epsilon_i\le -\gamma\}\nonumber\\
&=&\epsilon_i-\epsilon_iI\{|\epsilon_i|>\gamma\}+\gamma I\{\epsilon_i>\gamma\}-\gamma I\{\epsilon_i\le -\gamma\}\nonumber\\
&=&\epsilon_i-(\epsilon_i-\gamma)I\{\epsilon_i>\gamma\}-(\epsilon_i+\gamma)I\{\epsilon_i\le -\gamma\},
\ese
we have
\bse\label{eqn:1stmoment}
|E[H_\gamma'(\ep_i)|X_i]|&=&E[(\epsilon_i-\gamma)I\{\epsilon_i>\gamma\}+(-\epsilon_i-\gamma)I\{\epsilon_i\le -\gamma\}|X_i]\nonumber\\
&\le &E[|\epsilon_i|I\{|\epsilon_i|\ge \gamma\}|X_i]\nonumber\\
&\le& E[|\epsilon_i|^{1+\kappa}\gamma^{-\kappa}I\{|\epsilon|\ge \gamma\}|X_i]\nonumber\\
&\le &E[|\epsilon_i|^{1+\kappa}\gamma^{-\kappa}|X_i]\le v_\kappa \gamma^{-\kappa}.
\ese
Furthermore,
\bse\label{eqn:zeta}
  E[(H_\gamma'(\ep_i))^2|X_i]&=&E[\epsilon_i^2I\{|\epsilon_i|\le\gamma\}+\gamma^2 I\{|\epsilon_i|>\gamma\}|X_i]\nonumber\\
&\le& E[|\epsilon_i|^{1+\kappa}\gamma^{1-\kappa}|X_i]+\gamma^2 \cdot \frac{v_\kappa}{\gamma^{1+\kappa}}\nonumber\\
&\le &2v_\kappa\gamma^{1-\kappa}.
\ese
\hfill $\Box$

\section{Functional regression in an RKHS}\label{sec:rkhs}
We briefly discuss the case of kernel-based functional linear regression in the framework of reproducing kernel Hilbert space (RKHS), which has been adopted in \cite{caiyuan12,Zhang2020}. The readers are referred to \cite{aronszajn50,wahba90} for an introduction to the theory of RKHS. In this framework, we assume $f\in\calH\subset L^2(\calT)$ for some known  RKHS with its inner product denoted by $\langle\cdot,\cdot\rangle_\calH$, and estimate $f$ via
\be\label{eqn:modelrkhs}
\wh f=\argmin_{f\in\calH}\frac{1}{n}\sum_i H_\gamma(y_i-\langle X_i,f\rangle)+\lambda\|f\|_\calH^2,
\ee
where $\|.\|_\calH$ is the RKHS norm. It turns out the minimization problem above can be equivalently formulated in $L^2(\calT)$. Let $\calL_K$ be the operator induced by the kernel function $K$ associated with the RKHS, that is $\calL_K: g\in L^2(\calT)\rightarrow  \int K(.,t)g(t)dt,$ which is a positive-definite compact operator and we have $\calH=\calL_K^{1/2}(L^2(\calT))$ under the mild assumption that $K$ is positive-definite and $K\in L^2(\calT\times\calT)$ (Lemma 1.1.1 of \cite{wahba90}). Thus, there is a one-one mapping between $L^2(\calT)$ and $\calH$ given by $\calL_K^{1/2}: g\in L^2(\calT)\rightarrow f=\calL_K^{1/2}g\in\calH$, which is also an isometry: $\|g\|=\|f\|_\calH$. The minimization problem \eqref{eqn:modelrkhs} is equivalent to 
\bse
\wh g=\argmin_{g\in L^2(\calT)}\frac{1}{n}\sum_iH_\gamma(y_i-\langle \calL_K^{1/2}X_i,g\rangle)+\lambda\|g\|^2,
\ese
and then we can set $\wh f=\calL_K^{1/2}\wh g.$

With $\calL_K^{1/2}X_i$ as the functional predictor, and with $\Gamma_K=E[(\calL_K^{1/2}X_i)\otimes (\calL_K^{1/2}X_i)]=\calL_K^{1/2}\Gamma \calL_K^{1/2}$, our previous theoretical results yield the convergence rates of $\|(\Gamma_K+\lambda)^{1/2}(\wh g-g_0)\|$. We note that this error measure is actually equal to 
\bse
&&\langle(\calL_K^{1/2}\Gamma \calL_K^{1/2}+\lambda)(\wh g-g_0),\wh g-g_0\rangle\\
& =&\langle \Gamma \calL_K^{1/2}(\wh g-g_0),\calL_K^{1/2}(\wh g-g_0)\rangle +\lambda\|\wh g-g_0\|^2\\
&=&\|\Gamma^{1/2}(\wh f-f_0)\|^2 +\lambda\|f-f_0\|_\calH^2.
\ese
 In particular, $\|\Gamma^{1/2}(\wh f-f_0)\|$ is exactly the error measure for which the optimal statistical rate was established in \cite{caiyuan12}. Our rate is the same as that in \cite{caiyuan12} when $\kappa\ge 1$, although the latter paper only considers the case $r=0$, only used error measure $\|\Gamma^{1/2}(\wh f-f_0)\|$ (did not provide bounds for $\|\wh f-f_0\|_\calH$ as we can do here), and only considered convergence in probability.

\section{Conclusion}\label{sec:con}
This paper addresses the under-developed applications of robust estimation in functional linear regression, with a specific focus on Huber's loss incorporating a diverging robustness parameter. The proposed method, previously utilized in parametric models, stands out by aiming to robustly estimate the conditional mean, distinguishing itself from other robust methods like median regression that focus on estimating the conditional median, for example. The key advantage lies in the requirement of only a $(1+\kappa)$-th moment assumption ($\kappa>0$) on the noise distribution, with established error bounds matching the optimal rate in the least-squares scenario when $\kappa\ge 1$.

 By establishing convergence rates in probability for functional predictors with a finite 4-th moment, and providing finite-sample bounds with exponential tail for Gaussian functional predictors, the paper offers a comprehensive analysis of prediction error and $L^2$ error. These results also generalize to the realm of functional estimation within a reproducing kernel Hilbert space (RKHS).

\section*{Acknowledgements}  The research of Heng Lian is supported by NSFC 12371297 at CityU Shenzhen Research Institute, NSF of Jiangxi Province under Grant 20223BCJ25017, and by Hong Kong RGC general research fund 11300519, 11300721 and 11311822, and by CityU internal grant 7006014.

 
\appendix

\setcounter{equation}{0}
 
\renewcommand{\theequation}{S.\arabic{equation}}

  \bibliographystyle{jasa} 
  \bibliography{papers,books}

\end{document}